\documentclass[11pt]{article}
\usepackage{epic,latexsym,amssymb}
\usepackage{color}
\usepackage{tikz}

\newtheorem{theorem}{Theorem}
\newtheorem{lemma}{Lemma}
\newtheorem{quest}{Question}

\newtheorem{claim}{Claim}

\newtheorem{observation}{Observation}
\newtheorem{prop}{Proposition}

\newcommand{\QED}{$\Box$}
\newcommand{\smallqed}{{\tiny ($\Box$)}}

\newcommand{\modo}{{\rm mod}}
\newcommand{\cG}{\mathcal{G}}

\newcommand{\pn}{{\rm pn}}

\newcommand{\OOIR}{{\rm OOIR}}
\newcommand{\indM}{\alpha_{\rm ind}'}
\def\cp{\,\square\,}

\newcommand{\D}{{Dominator }}
\newcommand{\St}{{Staller }}

\newcommand{\epn}{{\rm epn}}
\newcommand{\ipn}{{\rm ipn}}

\newcommand{\gtd}{\gamma_{{\rm tg}}}

\newcommand{\grt}{\gamma_{\rm gr}^t}

\newcommand{\gti}{\gamma_t^{\rm i}} 
\newcommand{\gi}{\gamma_{\rm i}} 
\renewcommand{\gg}{\gamma_{\rm g}}

\newcommand{\proof}{\noindent\textbf{Proof. }}

\let\oldenumerate\enumerate
\renewcommand{\enumerate}{
  \oldenumerate
  \setlength{\itemsep}{0pt}
  \setlength{\parskip}{0pt}
  \setlength{\parsep}{0pt}
}

\begin{document}

\title{Indicated total domination game}

\author{$^1$Michael A. Henning\thanks{Research supported in part by the University of Johannesburg and the South African National Research Foundation} and $^2$Douglas F. Rall
\\ \\
$^1$Department of Mathematics and Applied Mathematics \\
University of Johannesburg \\
Auckland Park, 2006 South Africa\\
\small \tt Email: mahenning@uj.ac.za  \\
\\
$^2$Professor Emeritus of Mathematics \\
Furman University \\
Greenville, SC, USA\\
\small \tt Email: doug.rall@furman.edu}

\date{}
\maketitle

\begin{abstract}
A vertex $u$ in a graph $G$ totally dominates a vertex~$v$ if $u$ is adjacent to $v$ in $G$. A total dominating set of $G$ is a set $S$ of vertices of $G$ such that every vertex of $G$ is totally dominated by a vertex in $S$. The indicated total domination game is played on a graph $G$ by two players, Dominator and Staller, who take turns making a move. In each of his moves, Dominator indicates a vertex $v$ of the graph that has not been totally dominated in the previous moves, and Staller chooses (or selects) any vertex adjacent to $v$ that has not yet been played, and adds it to a set $D$ that is being built during the game. The game ends when every vertex is totally dominated, that is, when $D$ is a total dominating set of $G$. The goal of Dominator is to minimize the size of $D$, while Staller wants just the opposite. Providing that both players are playing optimally with respect to their goals, the size of the resulting set $D$ is the indicated total domination number of $G$, denoted by $\gti(G)$. In this paper we present several results on indicated total domination game. Among other results we prove that the indicated total domination number of a graph is bounded below by the well studied upper total domination number.
\end{abstract}


\newpage
\section{Introduction}

In 2010 Bre{\v{s}}ar, Klav{\v{z}}ar, and Rall~\cite{BrKlRa10} published the seminal paper on the domination game which belongs to the growing family of competitive optimization graph games. Domination games played on graphs are now very well studied in the literature. The subsequent rapid growth by the scientific community of research on domination games played on graphs inspired the recent book entitled ``Domination games played on graphs" by Bre\v{s}ar, Henning, Klav\v zar, and Rall~\cite{BrHeKlRa-21}, which presented the state of the art results at the time and shows that the area is rich for further research. In this paper, we study the total version of the indicated domination game.

A \emph{neighbor} of a vertex $v$ in $G$ is a vertex that is adjacent to $v$. The \emph{open neighborhood} of $v$ in $G$ is the set of neighbors of $v$, denoted $N_G(v)$. Thus, $N_G(v) = \{u \in V \colon uv \in E(G)\}$. The \emph{closed neighborhood of $v$} is the set $N_G[v] = \{v\} \cup N_G(v)$. For a set of vertices $S \subseteq V(G)$ and a vertex $v$ belonging to the set $S$, the \emph{open $S$-private neighborhood} is defined by $\pn(v,S) = \{w \in V \colon N_G(w) \cap S = \{v\}\}$. The \emph{$S$-external private neighborhood} of $v$ is the set $\epn(v,S) = \pn(v,S) \setminus S$, and the \emph{open $S$-internal private neighborhood} is the set $\ipn(v,S) = \pn(v,S) \cap S$. We note that $\pn(v,S) =  \epn(v,S) \cup \ipn(v,S)$. The \emph{closed neighborhood} of a set $S \subseteq V(G)$ is the set $\cup N_G[v]$ where the union is taken over all vertices $v \in S$.

 A vertex of degree~$1$ in a graph is called a \emph{leaf} and its neighbor a \emph{support vertex}.  An \emph{isolated vertex} in $G$ is a vertex of degree~$0$ in $G$. An \emph{isolate}-\emph{free graph} is a graph which contains no isolated vertex. A \emph{trivial graph} is the graph of order~$1$, and a \emph{nontrivial graph} has order at least~$2$. The \emph{join} of two graphs $G$ and $H$, denoted $G \oplus H$, is constructed from their disjoint union by adding edges making every vertex in $G$ adjacent to every vertex in $H$.  The subgraph of $G$ induced by a set $S \subseteq V(G)$ is denoted by $G[S]$. For an integer $k \ge 1$, we let $[k] = \{1,\ldots,k\}$ and $[k]_0 = \{0,1,\ldots,k\}$.

A vertex $u$ in a graph $G$ \emph{dominates} a vertex~$v$ if $u = v$ or $u$ is adjacent to $v$ in $G$. A \emph{dominating set} of $G$ is a set $S$ of vertices of $G$ such that every vertex in $G$ is dominated by a vertex in $S$. A vertex $u$ in a graph $G$ \emph{totally dominates} a vertex~$v$ if $u$ is adjacent to $v$ in $G$. A \emph{total dominating set}, abbreviated TD-set, of $G$ is a set $S$ of vertices of $G$ such that every vertex of $G$ is totally dominated by a vertex in $S$, that is, every vertex in $G$ has a neighbor in $S$. The \emph{total domination number} of $G$, denoted $\gamma_t(G)$, is the minimum cardinality of a TD-set in $G$, while the \emph{upper total domination number} of $G$, denoted $\Gamma_t(G)$, is the maximum cardinality of a minimal TD-set in $G$. A minimal TD-set of cardinality~$\Gamma_t(G)$ we call a $\Gamma_t$-\emph{set of $G$}. For other graph theory terminology not defined herein, the reader is referred to \cite{HaHeHe-23}, and for other recent books on domination in graphs, we refer the reader to~\cite{HaHeHe-20,HaHeHe-21,HeYe-book}.

\subsection{Domination games in graphs}

The \emph{domination game}, as introduced in~\cite{BrKlRa10}, is played on a graph $G$ by two players: Dominator and Staller. They alternate taking moves in which they select a vertex of $G$. A move is legal if the selected vertex dominates at least one vertex which is not already dominated by previously played vertices. The game ends when there are no legal moves, so when the set of played vertices is a dominating set of $G$. The goal of Dominator is to finish the game with minimum number of moves, while the aim of Staller is to maximize the number of moves. If both players play optimally, then the number of moves played on $G$ is the invariant (see~\cite{BrHeKlRa-21, fractional-2019}) the \emph{game domination number} $\gg(G)$, which is the number of moves on $G$ if Dominator starts the game. We refer to the domination game when \D makes the first move in the game as the $D$-game (standing for `Dominator start game').

The \emph{total domination game} was introduced in~\cite{Tdom-game-2015} as follows. Given a graph $G$, two players, called \emph{Dominator} and \emph{Staller}, take turns choosing a vertex from $G$. Each vertex chosen must totally dominate at least one vertex not totally dominated by the set of vertices previously chosen. Such a chosen vertex is called a \emph{legal move}. The game ends when there is no legal move available. Dominator wishes to minimize the number of vertices selected, while the goal of Staller is just the opposite. The \emph{game total domination number}, $\gtd(G)$, of $G$ is the number of moves played on $G$ when both players play optimally and Dominator starts the game. As before, we refer to the total domination game when \D makes the first move in the game as the $D$-game. The total domination game is now extensively studied in the literature; see, for example,~\cite{BrHeKlRa-21,Tdom-game-2017}.

A sequence of vertices in a graph $G$ is a \emph{total dominating sequence} if every vertex $v$ in the sequence totally dominates at least one vertex that was not totally dominated by any vertex that precedes $v$ in the sequence, and at the end all vertices of $G$ are totally dominated. While the length of a shortest such sequence is the total domination number $\gamma_t(G)$ of $G$, a total dominating sequence of maximum length is the \emph{Grundy total domination number}, $\grt(G)$, of~$G$.

\subsection{Open-open irredundant sets and induced matchings}

Fellows, Fricke, Hedetniemi, and Jacobs~\cite{FeFrHeJa-94} unified the study of private neighbors of vertices by forming a cohesive theory of private neighbors in graphs, yielding the so-called private neighbor cube. We mention here two such parameters defined in~\cite{FeFrHeJa-94}.

If each vertex in a set $S$ of vertices in a graph $G$ has a private neighbor inside $S$ other than itself, then the subgraph $G[S]$ induced by the set $S$ consists of disjoint copies of the complete graph $K_2$ on two vertices, that is, $G[S]$ is a $1$-regular subgraph of $G$. The set of edges in such an induced subgraph is called an \emph{induced matching} by Cameron~\cite{Ca-89} and Faudree, Gy\'{a}rf\'{a}s, Schelp, and Tuza~\cite{FaGyScTu-89} in 1989, and a \emph{strong matching} by Golumbic and Laskar~\cite{GoLa-93} in 1993. The maximum order of an induced matching in $G$ is the \emph{induced matching number}, denoted $\indM(G)$, of $G$. Induced matchings in graphs are now very well studied in the literature.

If each vertex in a set $S$ of vertices in a graph $G$ has a private neighbor other than itself, either inside $S$ or outside $S$, then such a set is coined an \emph{open}-\emph{open} \emph{irredundant set} in~\cite{FeFrHeJa-94}. Such a set satisfies the condition that for every vertex $v \in S$, the set $N(v) \setminus N(S \setminus \{v\})$ is not empty. The \emph{open}-\emph{open} \emph{irredundance number}, denoted $\OOIR(G)$, of $G$ is the maximum cardinality of an open-open irredundant set in $G$.

\subsection{Indicated total domination game}

Grzesik~\cite{gr-2012} proposed and first studied the \emph{indicated coloring game} which is a combinatorial game played on a graph $G$ by two players, and a fixed set $C$ of colors. In each round of the game Ann indicates an uncolored vertex, and Ben colors it using a color from $C$, obeying just the proper coloring rule. The goal of Ann is to achieve a proper coloring of the whole graph, while Ben is trying to prevent this. The minimum cardinality of the set of colors $C$ for which Ann has a winning strategy is called the \emph{indicated chromatic number}, $\chi_i(G)$, of a graph $G$.

Recently, Bre\v sar, Bujt\'as, Ir\v si\v c, Rall and Tuza~\cite{Ind-dom-game} introduced and studied the indicated domination game inspired by the indicated coloring game. The \emph{indicated domination game} is played on a graph $G$ by two players, \emph{Dominator} and \emph{Staller}, who take turns making a move. In each of his moves, Dominator indicates a vertex $v$ of the graph that has not been dominated in the previous moves, and Staller chooses (or selects) any vertex from the closed neighborhood of $v$ that has not yet been played, and adds it to a set $D$ that is being built during the game. The game ends when there is no undominated vertex left, that is, when $D$ is a dominating set. The goal of Dominator is to minimize the size of $D$, while Staller wants just the opposite. Providing that both players are playing optimally with respect to their goals, the size of the resulting set $D$ is the \emph{indicated domination number} of $G$, and is denoted by $\gi(G)$.

In this paper, we study the total version of the indicated domination game. The \emph{indicated total domination game} is played on an isolate-free graph $G$ by two players, \emph{Dominator} and \emph{Staller}, who take turns making a move. In each of his moves, Dominator indicates a vertex $v$ of the graph that has not been totally dominated in the previous moves, and Staller chooses (or selects) any vertex from the open neighborhood of $v$, and adds it to a set $D$ that is being built during the game. The game ends when every vertex is totally dominated, that is, when $D$ is a TD-set. The goal of Dominator is to minimize the size of $D$, while Staller wants just the opposite. Providing that both players are playing optimally with respect to their goals, the size of the resulting set $D$ is the \emph{indicated total domination number} of $G$, and is denoted by $\gti(G)$. Upon completion of the game, the set $D$ of vertices chosen by Staller is a TD-set of $G$ (possibly, $D = V(G)$). By considering the \emph{game tree} (which we do not define here) for the indicated total domination game, identical arguments as in~\cite[Chapter~1.2]{BrHeKlRa-21} show that the indicated total domination number is well-defined. Throughout this paper we only consider isolate-free graphs.

\section{The continuation principle}

A \emph{partially total dominated graph} is a graph together with a declaration that some vertices are already totally dominated. Since Dominator only indicates vertices of the graph that have not yet been totally dominated, he is not permitted to indicate any vertices that are declared already totally dominated. Adopting the notation in~\cite{Tdom-game-2017}, given a graph $G$ and a subset $S$ of vertices of $G$, we denote by $G|S$ the partially total dominated graph with $S$ as the set of declared vertices already totally dominated. Thus, $\gti(G|S)$ is the number of moves needed to finish the game on $G|S$ when Dominator and Staller both are playing optimally.

A key lemma in both the domination game and the total domination game (and other variants of certain domination games played on graphs) is the so-called \emph{Continuation Principle} first presented by Kinnersley, West, and Zamani in~\cite{bill-2013}. The Continuation Principle is a powerful tool for obtaining good upper bounds on game domination parameters.

\begin{figure}[ht!]
\begin{center}
\begin{tikzpicture}[scale=0.7,style=thick,x=0.9cm,y=0.9cm]
\def\vr{3pt} 
\path (-0.25,1) coordinate (u1);
\path (5.25,1) coordinate (u4);
\path (1,1) coordinate (v1);
\path (0.925,1) coordinate (v1p);
\path (2,0) coordinate (v2);
\path (3,0) coordinate (v3);
\path (4,1) coordinate (v4);
\path (3,2) coordinate (v5);
\path (2,2) coordinate (v6);
\draw (v1)--(v2)--(v3)--(v4)--(v5)--(v6)--(v1);
\draw (v1)--(u1);
\draw (v4)--(u4);
\draw (v1)--(v4);
\draw (v1) [fill=white] circle (\vr);
\draw (v2) [fill=white] circle (\vr);
\draw (v3) [fill=white] circle (\vr);
\draw (v4) [fill=white] circle (\vr);
\draw (v5) [fill=white] circle (\vr);
\draw (v6) [fill=white] circle (\vr);
\draw (u1) [fill=white] circle (\vr);
\draw (u4) [fill=white] circle (\vr);
\draw[anchor = north] (v2) node {{\small $v_2$}};
\draw[anchor = north] (v3) node {{\small $v_3$}};
\draw[anchor = south] (v5) node {{\small $v_5$}};
\draw[anchor = south] (v6) node {{\small $v_6$}};
\draw[anchor = south] (v1p) node {{\small $v_1$}};
\draw[anchor = south] (u1) node {{\small $u_1$}};
\draw[anchor = south] (v4) node {{\small $v_4$}};
\draw[anchor = south] (u4) node {{\small $u_4$}};
\end{tikzpicture}
\caption{A graph $G$}
\label{fig:counterexample}
\end{center}
\end{figure}
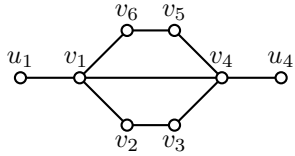

We show that the Continuation Principle unfortunately does not hold for the indicated total domination game. That is, we show that if $G$ is a graph and $A,B \subseteq V(G)$ with $B \subseteq A$, then it is not necessarily true that $\gti(G|A) \le \gti(G|B)$. Perhaps the simplest counterexample is to let $G$ be the graph obtained from a $6$-cycle $v_1v_2 \ldots v_6v_1$ by adding two new vertices $u_1$ and $u_4$, and adding the edges $u_1v_1$, $u_4v_4$ and $v_1v_4$. The graph $G$ is illustrated in Figure~\ref{fig:counterexample}.

\bigskip

In the graph $G$ shown in Figure~\ref{fig:counterexample}, let $B=\emptyset$ and let $A = \{u_1,u_4\}$. Suppose the indicated total domination game is played on $G|B$. We note that $G|B = G$, and so the game is played on the graph $G$.  In this game Dominator would first indicate vertex $u_1$, and Staller would be compelled to select vertex~$v_1$.  At this point the vertices in the set $\{u_1,v_2,v_4,v_6\}$ are totally dominated.  Dominator as his second move indicates vertex $u_4$ to which Staller must choose vertex $v_4$ and the game is complete. Thus, $\gti(G|B) \le 2$. Since $\gti(G|B) \ge \gamma_t(G) = 2$, we infer that $\gti(G|B) = 2$.

Now let the indicated total domination game be played on $G|A$.  If Dominator indicates vertex $v_1$, then Staller can choose vertex $u_1$. This leaves vertices $v_2, v_3,v_4,v_5$ and $v_6$ to be totally dominated, and that will require at least two more choices by Staller. If Dominator indicates vertex $v_4$, then by symmetry Staller can guarantee that at least three vertices are played (starting with the vertex $u_4$). If Dominator indicates vertex $v_2$, then Staller can choose vertex $v_1$.  This leaves vertices $v_1,v_3$ and $v_5$  to be totally dominated, and Staller can guarantee that at least two additional vertices are played. All other choices for Dominator to indicate follow from symmetry. Therefore, $\gti(G|A) \ge 3$. As observed earlier, $\gti(G|B) = 2$. Hence, $\gamma_t^i(G|B) < \gamma_t^i(G|A)$.

Since the \emph{Continuation Principle} does not hold for the indicated total domination game, this indicates that obtaining good upper bounds on the indicated total domination number is likely to be challenging.

\section{The upper total domination number}

In this section we show that the indicated total domination number of an isolate-free graph is at least the upper total domination of the graph. A fundamental property of minimal TD-sets was established by Cockayne, Dawes, and Hedetniemi~\cite{CoDaHe-80} in 1980.

\begin{lemma}{\rm (\cite{CoDaHe-80})}
\label{l:minimal-Tdom}
A {\rm TD}-set $S$ in a graph $G$ is a minimal {\rm TD}-set if and only if every vertex $v \in S$ has an open $S$-external private neighbor or an open $S$-internal private neighbor, that is, if and only if $|\epn(v,S)| \ge 1$ or $|\ipn(v,S)| \ge 1$.
\end{lemma}

As an application of Lemma~\ref{l:minimal-Tdom}, we can prove that the indicated total domination number of a graph is at least its upper total domination number.

\begin{prop}
\label{prop:upper-Tdom}
If $G$ is an isolate-free graph, then $\Gamma_t(G) \le \gti(G)$.
\end{prop}
\proof Let $S$ be a $\Gamma_t$-set of $G$, and so $S$ is a minimal TD-set of $G$ of cardinality $\Gamma_t(G)$. By Lemma~\ref{l:minimal-Tdom}, $|\epn(v,S)| \ge 1$ or $|\ipn(v,S)| \ge 1$ for every vertex $v \in S$, and so, $\pn(v,S) = \epn(v,S) \cup \ipn(v,S) \ne \emptyset$. Let $p = \Gamma_t(G)$ and define a partition $(V_1,V_2,\ldots,V_p)$ of the vertex set $V(G)$ into $p$ sets as follows. For $i \in [p]$, let $\pn(v_i,S) \subseteq V_i \subseteq N_G(v_i)$. Staller's strategy is to always select a vertex from $S$ to totally dominate the vertex indicated by Dominator. More precisely, whenever Dominator indicates a vertex $v$ to be totally dominated, Staller identifies the set $V_i$ that contains the vertex~$v$ for some $i \in [p]$ and selects the vertex $v_i \in S$. Since $\pn(w,S) \ne \emptyset$ for every vertex $w \in S$, this guarantees that upon completion of the game Staller selects all vertices in the set $S$, implying that $\gti(G) \ge |S| = \Gamma_t(G)$.~\QED

\medskip
The upper total domination number of a path $P_n$ of order~$n$ is established in~\cite{DoHeMc-07}.

\begin{prop}{\rm (\cite{DoHeMc-07})}
\label{prop:path}
For $n \ge 2$ an integer, $\Gamma_t(P_n) = 2 \lfloor \frac{n+1}{3} \rfloor$.
\end{prop}

We determine next the indicated total domination number of a path $P_n$ of order~$n$, and give a strategy for \D to play on a path. Recall that for $k \ge 1$, we let $[k] = \{1,\ldots,k\}$ and $[k]_0 = \{0\} \cup [k]$.

\begin{theorem}
\label{thm:path}
For $n \ge 2$, $\gti(P_n) = \Gamma_t(P_n)$.
\end{theorem}
\proof Let $T$ be a path $v_1v_2 \ldots v_n$ of order~$n \ge 2$. If $n = 2$, then it is immediate that $\gti(T) = \Gamma_t(T) = 2$. Hence we may assume that $n \ge 3$. By Proposition~\ref{prop:upper-Tdom}, $\gti(T) \ge \Gamma_t(T)$. Hence it suffices for us to prove that $\gti(T) \le \Gamma_t(T)$ from which we infer that $\gti(P_n) = \Gamma_t(P_n)$. We consider three cases.

\medskip
\emph{Case~1. $n \equiv 0 \, (\modo \, 3)$.} Thus, $n = 3k$ for some $k \ge 1$. Dominator's strategy is to indicate on his $(i+1)$st move the vertex $v_{3i+1}$ for $i \in [k-1]_0$. Thus on his first $k$ moves, \D indicates the vertices $v_1,v_4,\ldots,v_{3k-2}$ in turn. This forces Staller to play the vertex~$v_2$ on her first move, and if $k \ge 2$, then she is required to play either the vertex~$v_{3i}$ or $v_{3i+2}$ on her $(i+1)$st move for all $i \in [k-1]$. On each of her moves, Staller totally dominates two new vertices. Hence after the first $k$ moves of Staller, exactly $2k$ vertices on the path $T$ are totally dominated. For the remaining $k$ vertices on $T$ that are not yet totally dominated, \D simply indicates each such vertex on consecutive moves, thus guaranteeing that the game is complete after at most~$k$ additional moves of Staller. Thus, the game is finished in at most~$2k$ moves, and so $\gti(T) \le 2k = 2 \lfloor \frac{n+1}{3} \rfloor = \Gamma_t(P_n)$.

\medskip
\emph{Case~2. $n \equiv 1 \, (\modo \, 3)$.} Thus, $n = 3k+1$ for some $k \ge 1$. Dominator's strategy is to indicate on his $(i+1)$st move the vertex $v_{3i+1}$ for $i \in [k]_0$. Thus on his first $k+1$ moves, \D indicates the vertices $v_1,v_4,\ldots,v_{3k+1}$ in turn. This forces Staller to play the vertex~$v_2$ on her first move, the vertex~$v_{3i}$ or $v_{3i+2}$ on her $(i+1)$st move for all $i \in [k-1]$, and the vertex~$v_{3k}$ on her $(k+1)$st move. On each of her moves, Staller totally dominates two new vertices. Hence after the first $k+1$ moves of Staller, exactly $2(k+1)$ vertices on the path $T$ are totally dominated. For the remaining $k-1$ vertices on $T$ that are not yet totally dominated, \D indicates each such vertex on consecutive moves, thus guaranteeing that the game is complete after at most~$k-1$ additional moves of Staller. Thus, the game is finished in at most~$(k+1) + (k-1) = 2k$ moves, and so $\gti(T) \le 2k = 2 \lfloor \frac{n+1}{3} \rfloor = \Gamma_t(P_n)$.

\medskip
\emph{Case~3. $n \equiv 2 \, (\modo \, 3)$.} Thus, $n = 3k+2$ for some $k \ge 1$. Dominator's strategy is to indicate on his $i$th move the vertex $v_{3i-2}$ for $i \in [k]$. Thus on his first $k$ moves, \D indicates the vertices $v_1,v_4,\ldots,v_{3k-2}$ in turn. This forces Staller to play the vertex~$v_2$ on her first move, and to play the vertex~$v_{3i-3}$ or $v_{3i-1}$ on her $i$th move for $i \in [k] \setminus \{1\}$. On each of her moves, Staller totally dominates two new vertices. Hence after the first $k$ moves of Staller, exactly $2k$ vertices on the path $T$ are totally dominated. For the remaining $k+2$ vertices on $T$ that are not yet totally dominated, \D indicates each such vertex on consecutive moves, thus guaranteeing that the game is complete after at most~$k+2$ additional moves of Staller. Thus, the game is finished in at most~$k + (k+2) = 2k+2$ moves, and so $\gti(T) \le 2k + 2  = 2 \lfloor \frac{n+1}{3} \rfloor = \Gamma_t(P_n)$.~\QED

\medskip
A natural problem is to extend the path result in Theorem~\ref{thm:path} and to determine if the indicated total domination number is equal to the upper total domination in the class of trees. We state this formally as follows.

\begin{quest}
\label{quest:tree1}
Is it true that if $T$ is a nontrivial tree, then $\Gamma_t(T) = \gti(T)$?
\end{quest}

The following result gives a partial result in this direction.

\begin{prop}
\label{prop:tree1}
If $T$ is a nontrivial tree in which every vertex is a leaf or a support vertex, then $\gti(T) = \Gamma_t(T)$.
\end{prop}
\proof Let $T$ be a tree of order~$n \ge 2$ in which every vertex is a leaf or a support vertex. If $n = 2$, then the result is immediate. Hence we may assume that $n \ge 3$. If $T$ is a star $K_{1,n-1}$, then it is straightforward to check that $\Gamma_t(T) = \gti(T) = 2$. Hence, we may assume that $T$ has diameter at least~$3$, implying that $T$ contains at least two support vertices. Every TD-set of $T$ contains the set of support vertices of $T$ in order to totally dominate all leaves in $T$. Moreover, every minimal TD-set of $T$ is unique and consists of the support vertices of $T$, implying that $\Gamma_t(T) = s$ where here $s \ge 2$ denotes the number of support vertices in $T$. We show next that $\gti(T) = s$. Let $v_1,v_2,\ldots,v_s$ be the support vertices in $T$, and let $u_i$ be an arbitrary leaf neighbor of $v_i$ for $i \in [s]$. On Dominator's $i$th move, he indicates the leaf $u_i$ for $i \in [s]$. In order to totally dominate the leaf $u_i$, Staller is required to play the vertex $v_i$ on her $i$th move. The resulting set of $s$ vertices played by Staller is the set of support vertices in $T$, which is a TD-set of $T$ and completes the game. Hence, Dominator can guarantee that the indicated total domination game finishes in at most~$s$ moves, that is, $\gti(T) \le s = \Gamma_t(T)$. By Proposition~\ref{prop:upper-Tdom}, $\Gamma_t(G) \le \gti(G)$. Consequently, $\Gamma_t(G) = \gti(G)$.~\QED

\medskip
We show next that there exists connected graphs $G$ satisfying $\Gamma_t(G) < \gamma_t^i(G)$. For this purpose, for a given graph $G$ and an integer $k \ge 2$, the \emph{$k^{\rm th}$-power graph} of $G$, denoted $G^k$, is the graph with the same vertex set as $G$ and where two vertices $u$ and $v$ are adjacent in $G^k$ if $d_G(u,v) \le k$. We show in the following result that the $k$th power of a cycle of order~$2k+3$ has upper total domination number less than its indicated total domination number.

\begin{prop}
\label{prop:powerk}
If $k \ge 2$ and $G=C_{2k+3}^k$, then $\Gamma_t(G) < \gamma_t^i(G)$.
\end{prop}
\proof
Let $C$ be the cycle $C_{2k+3}$ with $V(C) =\{v_1,\ldots,v_{2k+3}\}$ and $E(C) = \{ v_iv_{i+1} \colon i \in [2k+3] \}$ where the indices are computed modulo $2k+3$, and consider the $k$th power $G = C^k$ of the cycle $C$. First we claim that no minimal TD-set of $G$ contains two consecutive vertices from the original cycle $C$.  Without loss of generality suppose, to the contrary, that $\{v_1,v_2\}$ is a subset of a minimal TD-set $A$ of $G$. We note that $N_G(\{v_1,v_2\}) = V(G) \setminus \{v_{k+3}\}$. This implies that $A$ contains a neighbor, $v_j$, of $v_{k+3}$ in $G$. Thus, $\{v_1,v_2,v_j\} \subseteq A$, where $3 \le j \le 2k+3$ and $j \ne k+3$. If $j = k+4$, then since $d_C(v_1,v_{k+4}) = k$ and $d_C(v_2,v_{k+4}) = k+1$, we infer that $N(v_2) \setminus N(A \setminus \{v_2\}) = \emptyset$, and so $\epn(v_2,A) = \ipn(v_2,A) = \emptyset$. If $j = k+2$, then since $d_C(v_1,v_{k+2}) = k+1$ and $d_C(v_2,v_{k+2}) = k$, we infer that $N(v_1) \setminus N(A \setminus \{v_1\}) = \emptyset$, and so $\epn(v_1,A) = \ipn(v_1,A) = \emptyset$. If $k+5 \le j \le 2k+3$, then $N(v_1) \setminus N(A \setminus \{v_1\})=\emptyset$, and so $\epn(v_1,A) = \ipn(v_1,A) = \emptyset$. If $3 \le j \le k+1$, then $N(v_2) \setminus N(A \setminus \{v_2\}) = \emptyset$, and so $\epn(v_2,A) = \ipn(v_2,A) = \emptyset$. In all cases, we contradict the minimality of the TD-set $A$ as given by Lemma~\ref{l:minimal-Tdom}. Hence, no minimal TD-set of $G$ contains two consecutive vertices from the original cycle $C$. However if $1 \le i < j \le 2k+3$ and $k \ge d_{C}(v_i,v_j) > 1$, then $\{v_i,v_j\}$ is a minimal TD-set of $G = C^k$ and, in fact, every minimal TD-set of $G$ is formed this way.  We infer that $\Gamma_t(G)=2$.

When the indicated total domination game is played on $G$ we may assume by symmetry that $v_1$ is the first vertex chosen by Staller.  The vertices $v_1,v_{k+2}$ and $v_{k+3}$ are the only three vertices in $G$ not totally dominated by the vertex~$v_1$.  If Dominator now indicates $v_1$ or $v_{k+2}$ on his second move, then Staller can select $v_2$ which leaves $v_{k+3}$ not totally dominated.  On the other hand, if Dominator now indicates $v_{k+3}$ on his second move, then Staller can choose $v_{2k+3}$, which leaves $v_{k+2}$ not totally dominated.  Therefore, Staller can force at least three vertices to be played upon completion of the game, implying that $\gamma_t^i(G)=3$.~\QED

\medskip
By Proposition~\ref{prop:powerk}, there exists connected graphs $G$ of arbitrarily large minimum degree satisfying $\Gamma_t(G) < \gamma_t^i(G)$. We show in the next section that there exist connected graphs $G$ such that $\gamma_t^i(G)$ can be arbitrarily larger than $\Gamma_t(G)$.

\section{The open-open irredundance number}
\label{S:ooir}

Let $G$ be an isolate-free graph. From the definition of an open-open irredundant set, if $S$ is a minimal TD-set in $G$, then $S$ is an open-open irredundant set. In particular, every $\Gamma_t$-set of $G$ is an open-open irredundant set of $G$. Hence as observed in~\cite{FeFrHeJa-94}, it holds that $\Gamma_t(G) \le \OOIR(G)$. By Proposition~\ref{prop:upper-Tdom}, we have $\Gamma_t(G) \le \gti(G)$. Hence it is a natural question to ask whether $\gti(G) \le \OOIR(G)$. We show in this section that in general there is no relation between the indicated total domination number, $\gti(G)$, of a graph $G$ and the open-open irredundance number, $\OOIR(G)$, of $G$.

Let $\cG$ be the family of graphs constructed as follows. Let $H$ be the graph obtained from the join of two vertex disjoint copies $Q_1 \colon u_1v_1w_1x_1$ and $Q_2 \colon u_2v_2w_2x_2$ of a path $P_4$ of order~$4$, and so $H = P_4 \oplus P_4$.  For $k \ge 1$, let $H_1, \ldots, H_k$ be $k$ vertex disjoint copies of $H$. Further, let $Q_{i,1} \colon u_{i,1}v_{i,1}w_{i,1}x_{i,1}$ and $Q_{i,2} \colon u_{i,2}v_{i,2}w_{i,2}x_{i,2}$ be the paths in $H_i$ corresponding to the paths $Q_1$ and $Q_2$ in $H$ for $i \in [k]$. Let $G_1 = H_1$, and for $k \ge 2$ let $G_k$ be obtained from the disjoint union of the graphs $H_1, \ldots, H_k$ by adding the edges $w_{i,2}w_{i+1,1}$ for all $i \in [k]$ where addition is taken modulo~$k$. The graph $G_k$ is illustrated in Figure~\ref{fig:join}. Let $\cG = \{G_k \colon k \ge 1\}$.

\bigskip
\begin{figure}[ht!]
\begin{center}
\begin{tikzpicture}[scale=0.7,style=thick,x=0.9cm,y=0.9cm]
\def\vr{3pt} 
\path (0,-2) coordinate (u1);
\path (0,0) coordinate (v1);
\path (0,2) coordinate (w1);
\path (0,4) coordinate (x1);
\path (0,5) coordinate (x1p);
\path (2,-2) coordinate (u2);
\path (2,0) coordinate (v2);
\path (2,2) coordinate (w2);
\path (2,4) coordinate (x2);

\path (5,-2) coordinate (u3);
\path (5,0) coordinate (v3);
\path (5,2) coordinate (w3);
\path (5,4) coordinate (x3);
\path (7,-2) coordinate (u4);
\path (7,0) coordinate (v4);
\path (7,2) coordinate (w4);
\path (7,4) coordinate (x4);

\path (10,-2) coordinate (u5);
\path (10,0) coordinate (v5);
\path (10,2) coordinate (w5);
\path (10,4) coordinate (x5);
\path (12,-2) coordinate (u6);
\path (12,0) coordinate (v6);
\path (12,2) coordinate (w6);
\path (12,4) coordinate (x6);

\path (17,-2) coordinate (u7);
\path (17,0) coordinate (v7);
\path (17,2) coordinate (w7);
\path (17,4) coordinate (x7);
\path (19,-2) coordinate (u8);
\path (19,0) coordinate (v8);
\path (19,2) coordinate (w8);
\path (19,4) coordinate (x8);
\path (19,5) coordinate (x8p);
%
\draw (w2)--(w3); \draw (w4)--(w5); \draw(w6)-- (12.5,2); \draw (16.5,2)--(w7);
\foreach \i in {1,2,3,4,5,6,7,8}
{\draw (x\i)--(w\i); \draw (w\i)--(v\i); \draw (v\i)--(u\i);}
\draw (x1)--(x2); \draw (x1)--(w2); \draw (x1)--(v2); \draw (x1)--(u2); \draw (w1)--(x2); \draw (w1)--(w2); \draw (w1)--(v2); \draw (w1)--(u2);
\draw (v1)--(x2); \draw (v1)--(w2); \draw (v1)--(v2); \draw (v1)--(u2); \draw (u1)--(x2); \draw (u1)--(w2); \draw (u1)--(v2); \draw (u1)--(u2);
\draw (x3)--(x4); \draw (x3)--(w4); \draw (x3)--(v4); \draw (x3)--(u4); \draw (w3)--(x4); \draw (w3)--(w4); \draw (w3)--(v4); \draw (w3)--(u4);
\draw (v3)--(x4); \draw (v3)--(w4); \draw (v3)--(v4); \draw (v3)--(u4); \draw (u3)--(x4); \draw (u3)--(w4); \draw (u3)--(v4); \draw (u3)--(u4);
\draw (x5)--(x6); \draw (x5)--(w6); \draw (x5)--(v6); \draw (x5)--(u6); \draw (w5)--(x6); \draw (w5)--(w6); \draw (w5)--(v6); \draw (w5)--(u6);
\draw (v5)--(x6); \draw (v5)--(w6); \draw (v5)--(v6); \draw (v5)--(u6); \draw (u5)--(x6); \draw (u5)--(w6); \draw (u5)--(v6); \draw (u5)--(u6);
\draw (x7)--(x8); \draw (x7)--(w8); \draw (x7)--(v8); \draw (x7)--(u8); \draw (w7)--(x8); \draw (w7)--(w8); \draw (w7)--(v8); \draw (w7)--(u8);
\draw (v7)--(x8); \draw (v7)--(w8); \draw (v7)--(v8); \draw (v7)--(u8); \draw (u7)--(x8); \draw (u7)--(w8); \draw (u7)--(v8); \draw (u7)--(u8);
\draw (w1) to[out=180,in=180, distance=1.75cm] (x1p);
\draw (x1p)--(x8p);
\draw (w8) to[out=0,in=0, distance=1.75cm] (x8p);
\foreach \i in {1,2,3,4,5,6,7,8}
{  \draw (u\i)  [fill=white] circle (\vr);
\draw (v\i)  [fill=white] circle (\vr);
\draw (w\i)  [fill=white] circle (\vr);
\draw (x\i)  [fill=white] circle (\vr);}

\draw (13.75,2) [fill=black] circle (1.75pt);
\draw (14.5,2) [fill=black] circle (1.75pt);
\draw (15.25,2) [fill=black] circle (1.75pt);
\draw[anchor = east] (x1) node {{\small $x_{1,1}$}};
\draw[anchor = east] (x3) node {{\small $x_{2,1}$}};
\draw[anchor = east] (x5) node {{\small $x_{3,1}$}};
\draw[anchor = east] (x7) node {{\small $x_{k,1}$}};
\draw[anchor = west] (x2) node {{\small $x_{1,2}$}};
\draw[anchor = west] (x4) node {{\small $x_{2,2}$}};
\draw[anchor = west] (x6) node {{\small $x_{3,2}$}};
\draw[anchor = west] (x8) node {{\small $x_{k,2}$}};
\draw (w1)+(-.5,-.3) node {{\small $w_{1,1}$}};
\draw (w3)+(-.5,-.3) node {{\small $w_{2,1}$}};
\draw (w5)+(-.5,-.3) node {{\small $w_{3,1}$}};
\draw (w7)+(-.7,-.3) node {{\small $w_{k,1}$}};
\draw (w2)+(.7,-.3) node {{\small $w_{1,2}$}};
\draw (w4)+(.7,-.3) node {{\small $w_{2,2}$}};
\draw (w6)+(.7,-.3) node {{\small $w_{3,2}$}};
\draw (w8)+(.7,-.3) node {{\small $w_{k,2}$}};
\draw[anchor = east] (v1) node {{\small $v_{1,1}$}};
\draw[anchor = east] (v3) node {{\small $v_{2,1}$}};
\draw[anchor = east] (v5) node {{\small $v_{3,1}$}};
\draw[anchor = east] (v7) node {{\small $v_{k,1}$}};
\draw[anchor = west] (v2) node {{\small $v_{1,2}$}};
\draw[anchor = west] (v4) node {{\small $v_{2,2}$}};
\draw[anchor = west] (v6) node {{\small $v_{3,2}$}};
\draw[anchor = west] (v8) node {{\small $v_{k,2}$}};

\draw[anchor = east] (u1) node {{\small $u_{1,1}$}};
\draw[anchor = east] (u3) node {{\small $u_{2,1}$}};
\draw[anchor = east] (u5) node {{\small $u_{3,1}$}};
\draw[anchor = east] (u7) node {{\small $u_{k,1}$}};
\draw[anchor = west] (u2) node {{\small $u_{1,2}$}};
\draw[anchor = west] (u4) node {{\small $u_{2,2}$}};
\draw[anchor = west] (u6) node {{\small $u_{3,2}$}};
\draw[anchor = west] (u8) node {{\small $u_{k,2}$}};
\end{tikzpicture}
\caption{The graph $G_k$ in the family $\cG$}
\label{fig:join}
\end{center}
\end{figure}

\begin{prop}
\label{prop:join}
For $k \ge 1$, we have $\gamma_t^i(G_k)=3k$, $\Gamma_t(G_k)=2k$ and $\OOIR(G_k)=2k$.
\end{prop}
\proof
For $k \ge 1$, consider the graph $G = G_k \in \cG$. Let the indicated total domination game be played on $G$.  We provide a strategy for Staller that will ensure at least $3k$ vertices are chosen. By symmetry we may assume that in his first move Dominator indicates a vertex that belongs to the path $Q_{1,1}$. Staller responds by playing $x_{1,2}$, which totally dominates all vertices on the path $Q_{1,1}$ and the vertex~$w_{1,2}$.  Furthermore, later in the game Staller will only play vertices from the path $Q_{1,2}$ when Dominator indicates a vertex in $Q_{1,2}$ that is not yet totally dominated. This implies that at least three vertices from $Q_{1,2}$ will be played by Staller since both $v_{1,2}$ and $w_{1,2}$ are support vertices in the path $Q_{1,2}$.

As play progresses and Dominator first indicates a vertex in $H_j$ for some $j \in [k] \setminus \{1\}$, Staller will choose a vertex in a manner similar to how she responded to Dominator's first play.  That is, if Dominator indicates a vertex that belongs to the path $Q_{j,1}$, then Staller plays the vertex $x_{j,2}$. In this case when Dominator indicates any further vertices in $Q_{j,2}$, Staller will choose a vertex that belongs to the path $Q_{j,2}$. On the other hand, if the first vertex Dominator indicates from $H_j$ is a vertex that belongs to the path $Q_{j,2}$, then Staller chooses the vertex $x_{j,1}$ and ensures that two more additional vertices from $Q_{j,1}$ are chosen in the remainder of the game. This strategy ensures that at least three vertices will be chosen by Staller in each of the $k$ copies of $H$ in the graph $G_k$. Dominator can prevent more than three vertices being chosen from $H_j$, for each $j \in[k]$, by playing as follows.  The first vertex he indicates in $H_j$ should be $x_{j,1}$.  A short analysis shows that it is then to Staller's advantage to select either $x_{j,2}$ or $u_{j,2}$. By indicating $x_{j,2}$ and then $u_{j,2}$ on his next two moves Dominator ensures that at most three vertices will be selected from $H_j$. We therefore infer that Staller has a strategy to play exactly three vertices from every copy of $H = P_4 \oplus P_4$ in $G_k$, implying that $\gamma_t^i(G_k)=3k$.

Let $D$ be a $\Gamma_t$-set of $G_k$, and so $D$ is a minimal TD-set of $G_k$ of maximum cardinality. Let $D_i = D \cap V(H_i)$ for all $i \in [k]$. We show that $|D_i| \le 2$. If $D_i$ contains a vertex from $Q_{i,1}$ and a vertex from $Q_{i,2}$, then two such vertices form a set that totally dominates all vertices of $H_i$, and by the minimality of the set $D$ we infer that $|D_i| = 2$. If $D_i$ contains vertices from exactly one of $Q_{i,1}$ and $Q_{i,2}$, say from $Q_{i,1}$, then in order to totally dominate the vertices $u_{i,1}$ and $x_{i,1}$ the set $D_i$ contains the vertices $v_{i,1}$ and $w_{i,1}$, respectively. However, the set $\{v_{i,1},w_{i,1}\}$ totally dominates all vertices of $H_i$, and once again by the minimality of the set $D$ we infer that $|D_i| = 2$. This is true for all $i \in [k]$, and so $\Gamma_t(G_k) \le 2k$. Since the set $\{x_{1,1},x_{1,2}, \ldots, x_{k,1},x_{k,2}\}$, for example, is a minimal TD-set of $G_k$, we note that $\Gamma_t(G_k) \ge 2k$. Consequently, $\Gamma_t(G_k) = 2k$. Similarly, any open-open irredundant set in $G_k$ can contain at most two vertices from $H_i$ for each $i \in [k]$, and so $\OOIR(G_k) \le 2k$. As observed earlier, $2k = \Gamma_t(G) \le \OOIR(G) \le 2k$. Consequently, $\OOIR(G) = 2k$.~\QED

\medskip
By Proposition~\ref{prop:join}, there exist connected graphs $G$ such that $\gamma_t^i(G)$ can be arbitrarily larger than $\Gamma_t(G)$ and $\OOIR(G)$.
We show next that there exist connected graphs $G$ such that $\OOIR(G)$ can be arbitrarily larger than $\gamma_t^i(G)$. Let $k$ be a positive integer at least~$5$, and let $F_{k,1}$ and $F_{k,2}$ be two disjoint copies of the complete graph $K_k$, where $V(F_{k,1}) = \{u_1,u_2,\ldots, u_k\}$ and $V(F_{k,2}) = \{v_1,v_2,\ldots, v_k\}$. Let $F_k$ be obtained from the disjoint union of $F_{k,1}$ and $F_{k,2}$ by adding the $k-1$ edges $u_iv_i$ for $i \in [k-1]$. Thus, $F_k$ is obtained from the prism $K_k \cp K_2$ of a complete graph $K_k$ by removing one of the added edges between a pair of corresponding vertices in the copies of the complete graph, that is, $F_k \cong (K_k \cp K_2)-u_kv_k$. The graph $F_k$ is illustrated in Figure~\ref{fig:CartesianProduct}, where for clarity we omit the edges in the complete graphs $F_{k,1}$ and $F_{k,2}$.

\begin{figure}[ht!]
\begin{center}
\begin{tikzpicture}[scale=0.7,style=thick,x=0.9cm,y=0.9cm]
\def\vr{3pt} 
\path (1,1) coordinate (u1);
\path (1,0.5) coordinate (u11);
\path (1,2) coordinate (u12);
\path (2.5,1) coordinate (u2);
\path (4,1) coordinate (u3);
\path (8,1) coordinate (u4);
\path (9.5,1) coordinate (u5);
\path (9.5,0.5) coordinate (u51);
\path (9.5,2) coordinate (u52);
\path (1,-1) coordinate (v1);
\path (1,-0.5) coordinate (v11);
\path (1,-2) coordinate (v12);
\path (2.5,-1) coordinate (v2);
\path (4,-1) coordinate (v3);
\path (8,-1) coordinate (v4);
\path (9.5,-1) coordinate (v5);
\path (9.5,-0.5) coordinate (v51);
\path (9.5,-2) coordinate (v52);


\foreach \i in {1,2,3,4}
{ \draw (u\i)--(v\i);
}

\foreach \i in {1,2,3,4,5}
{  \draw (u\i)  [fill=white] circle (\vr); \draw (v\i)  [fill=white] circle (\vr);
 }

\draw (5.25,1) [fill=black] circle (1.5pt);
\draw (6,1) [fill=black] circle (1.5pt);
\draw (6.75,1) [fill=black] circle (1.5pt);
\draw (5.25,-1) [fill=black] circle (1.5pt);
\draw (6,-1) [fill=black] circle (1.5pt);
\draw (6.75,-1) [fill=black] circle (1.5pt);

\draw (u11) to[out=180,in=180, distance=1.65cm] (u12);
\draw (v11) to[out=180,in=180, distance=1.65cm] (v12);
\draw (u51) to[out=0,in=0, distance=1.65cm] (u52);
\draw (v51) to[out=0,in=0, distance=1.65cm] (v52);
\draw (u11)--(u51);
\draw (u12)--(u52);
\draw (v11)--(v51);
\draw (v12)--(v52);

\draw[anchor = south] (u1) node {{\small $u_1$}};
\draw[anchor = south] (u2) node {{\small $u_2$}};
\draw[anchor = south] (u3) node {{\small $u_3$}};
\draw[anchor = south] (u4) node {{\small $u_{k-1}$}};
\draw[anchor = south] (u5) node {{\small $u_k$}};
\draw[anchor = north] (v1) node {{\small $v_1$}};
\draw[anchor = north] (v2) node {{\small $v_2$}};
\draw[anchor = north] (v3) node {{\small $v_3$}};
\draw[anchor = north] (v4) node {{\small $v_{k-1}$}};
\draw[anchor = north] (v5) node {{\small $v_k$}};

\draw (6,1.5) node {{\small $K_k$}};
\draw (6,-1.5) node {{\small $K_k$}};

\end{tikzpicture}
\caption{The graph $F_k = (K_k \cp K_2) - u_kv_k$}
\label{fig:CartesianProduct}
\end{center}
\end{figure}

\newpage
\begin{prop}
\label{prop:cartesian}
For $k \ge 5$, we have $\gamma_t^i(F_k)=4$ and $\OOIR(F_k)=k-1$.
\end{prop}
\proof The set $\{u_1,u_2,\ldots, u_{k-1}\}$ is an open-open irredundant set of $F_k$, which implies that $\OOIR(F_n) \ge k-1$.  Suppose $A$ is an open-open irredundant subset of $V(F_n)$ such that $|A| \ge k$. If $|A \cap \{u_1,u_2,\ldots, u_k\}| \ge 1$ and $|A \cap \{v_1,v_2,\ldots, v_k\}| \ge 1$, then neither of these intersections has cardinality more than~$2$ and thus $|A|\le 4$, which is a contradiction.  We may thus assume that $A = \{u_1,u_2,\ldots, u_k\}$.  This is also a contradiction since $N(u_k) \setminus N(A \setminus \{u_k\}) = \emptyset$. Therefore, $\OOIR(F_k)=k-1$.

We show next that $\gamma_t^i(F_k)=4$. Let the indicated total domination game be played on the graph $F_k$. Suppose that Dominator indicates $u_k$ on his first move. Staller can only choose a vertex from $\{u_1,u_2,\ldots, u_{k-1}\}$ to totally dominate $u_k$.  Without loss of generality we may assume Staller plays $u_1$, and this vertex totally dominates $\{v_1,u_2,u_3,\ldots,u_k\}$. Suppose that Dominator indicates $v_k$ on his second move.  Staller could choose any vertex from $\{v_1,v_2,\ldots, v_{k-1}\}$ to totally dominate $v_k$.  If she chooses $v_1$, then the game ends. Any other choice by Staller, say $v_j$ where $2 \le j \le k-1$, leaves just two vertices, namely, $u_1$ and $v_j$, not totally dominated.  Thus, at most two more vertices can be chosen by Staller.  Regardless of which of these two vertices are indicated by Dominator on his third move, Staller can choose the appropriate vertex from $\{u_k,v_k\}$ and thus make sure the game ends with four vertices being chosen. Suppose next that Dominator indicates $v_j$ on his second move, where $2 \le j \le k-1$. In this case, Staller chooses the vertex $v_k$, thereby totally dominating all vertices except for $u_1$ and $v_k$. If Dominator indicates $u_1$, then Staller chooses the vertex $u_k$, while if Dominator indicates $v_k$, then Staller chooses the vertex $v_2$. Regardless of what vertex \D indicates, Staller can force four vertices to be chosen before the game ends.

Suppose that \D indicates a vertex different from $u_k$ and $v_k$. By symmetry and renaming vertices if necessary, we may assume that \D indicates vertex $u_1$. In this case, Staller chooses the vertex $u_k$, thereby totally dominating vertices in the set $\{u_1,u_2,\ldots, u_{k-1}\}$. Dominator must then indicate a vertex from the set $\{v_1,v_2,\ldots,v_k,u_k\}$ on his second move. However, regardless of what vertex \D indicates, Staller can force four vertices to be chosen before the game ends. Therefore, $\gamma_t^i(F_k)= 4$.~\QED

\medskip
By Proposition~\ref{prop:cartesian}, there exist connected graphs $G$ such that $\OOIR(G)$ can be arbitrarily larger than $\gamma_t^i(G)$.

\section{The induced matching number}
\label{S:ind-math}

The following relation between the open-open irredundance number and the induced matching number is given in~\cite{FeFrHeJa-94}.

\begin{theorem}{\rm (\cite{FeFrHeJa-94})}
\label{thm:cube-paper}
If $G$ is an isolate-free graph, then $2\indM(G) \le \OOIR(G)$. Moreover if $G$ is a bipartite graph, then $2\indM(G) = \OOIR(G)$.
\end{theorem}

As observed earlier, $\Gamma_t(G) \le \OOIR(G)$ holds for all isolate-free graphs $G$. Hence as a consequence of Theorem~\ref{thm:cube-paper}, if $G$ is an isolate-free bipartite graph, then $\Gamma_t(G) \le 2\indM(G)$. Recall that by Proposition~\ref{prop:upper-Tdom}, if $G$ is an isolate-free graph, then $\Gamma_t(G) \le \gti(G)$. A natural problem is to determine if the parameters $\gti(G)$ and $2\indM(G)$ are related.

Let $\cG = \{G_k \colon k \ge 1\}$ be the family of graphs constructed earlier (see Figure~\ref{fig:join} for an illustration of the graph $G_k$ that belongs to the family~$\cG$). For $k \ge 1$ if $G_k \in \cG$, then  by Proposition~\ref{prop:join} we have $\gamma_t^i(G_k)=3k$. Since each of the $k$ copies of the graph $H = P_4 \oplus P_4$ used to build the graph $G_k$ is dense and has induced matching number equal to~$1$, we observe that $\indM(G_k) = k$. This yields the following result.

\begin{prop}
\label{prop:ind-match1}
For $k \ge 1$, we have $\gamma_t^i(G_k)=3k$ and $2\indM(G_k) = 2k$.
\end{prop}

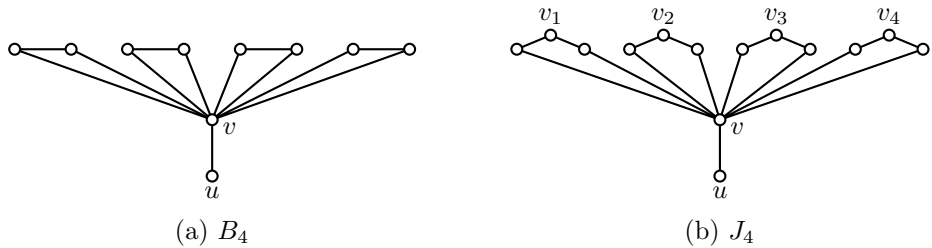
\begin{figure}[htb]
\begin{center}
\begin{tikzpicture}[scale=.75,style=thick,x=1cm,y=1cm]
\def\vr{2.75pt}
\path (0,2.25) coordinate (v1);
\path (1,2.25) coordinate (v2);
\path (2,2.25) coordinate (v3);
\path (3,2.25) coordinate (v4);
\path (4,2.25) coordinate (v5);
\path (5,2.25) coordinate (v6);
\path (6,2.25) coordinate (v7);
\path (7,2.25) coordinate (v8);
\path (3.5,1) coordinate (v);
\path (3.5,0.85) coordinate (vp);
\path (3.5,0) coordinate (u);
\draw (v)--(v1)--(v2)--(v)--(u);
\draw (v)--(v3)--(v4)--(v);
\draw (v)--(v5)--(v6)--(v);
\draw (v)--(v7)--(v8)--(v);
%
\draw (v1) [fill=white] circle (\vr);
\draw (v2) [fill=white] circle (\vr);
\draw (v3) [fill=white] circle (\vr);
\draw (v4) [fill=white] circle (\vr);
\draw (v5) [fill=white] circle (\vr);
\draw (v6) [fill=white] circle (\vr);
\draw (v7) [fill=white] circle (\vr);
\draw (v8) [fill=white] circle (\vr);
\draw (u) [fill=white] circle (\vr);
\draw (v) [fill=white] circle (\vr);
\draw[anchor = north] (u) node {{\small $u$}};
\draw[anchor = west] (vp) node {{\small $v$}};
\path (8.9,2.25) coordinate (v1);
\path (9.5,2.5) coordinate (v12);
\path (10.1,2.25) coordinate (v2);
\path (10.9,2.25) coordinate (v3);
\path (11.5,2.5) coordinate (v34);
\path (12.1,2.25) coordinate (v4);
\path (12.9,2.25) coordinate (v5);
\path (13.5,2.5) coordinate (v56);
\path (14.1,2.25) coordinate (v6);
\path (14.9,2.25) coordinate (v7);
\path (15.5,2.5) coordinate (v78);
\path (16.1,2.25) coordinate (v8);
\path (12.5,1) coordinate (v);
\path (12.5,0.85) coordinate (vp);
\path (12.5,0) coordinate (u);
\draw (v)--(v1)--(v12)--(v2)--(v)--(u);
\draw (v)--(v3)--(v34)--(v4)--(v);
\draw (v)--(v5)--(v56)--(v6)--(v);
\draw (v)--(v7)--(v78)--(v8)--(v);
%
\draw (v1) [fill=white] circle (\vr);
\draw (v2) [fill=white] circle (\vr);
\draw (v3) [fill=white] circle (\vr);
\draw (v4) [fill=white] circle (\vr);
\draw (v5) [fill=white] circle (\vr);
\draw (v6) [fill=white] circle (\vr);
\draw (v7) [fill=white] circle (\vr);
\draw (v8) [fill=white] circle (\vr);
\draw (u) [fill=white] circle (\vr);
\draw (v) [fill=white] circle (\vr);
\draw (v12) [fill=white] circle (\vr);
\draw (v34) [fill=white] circle (\vr);
\draw (v56) [fill=white] circle (\vr);
\draw (v78) [fill=white] circle (\vr);
\draw[anchor = north] (u) node {{\small $u$}};
\draw[anchor = west] (vp) node {{\small $v$}};
\draw[anchor = south] (v12) node {{\small $v_1$}};
\draw[anchor = south] (v34) node {{\small $v_2$}};
\draw[anchor = south] (v56) node {{\small $v_3$}};
\draw[anchor = south] (v78) node {{\small $v_4$}};
\draw (3.5,-1) node {{\small (a) $B_4$}};
\draw (12.5,-1) node {{\small (b) $J_4$}};
\end{tikzpicture}
\caption{The graphs $B_4$ and $J_4$}
\label{fig:B4}
\end{center}
\end{figure}

Moreover for $k \ge 1$, let $B_k$ be the graph obtained from $k$ vertex disjoint copies of $K_3$ by selecting one vertex from each copy of $K_3$ and identifying these vertices into one new vertex~$v$, and then adding a new vertex $u$ and adding the edge $uv$. The graph $B_4$, for example, is illustrated in Figure~\ref{fig:B4}(a). Dominator can force the game to be completed after two moves by indicating the vertex~$u$ on his first move, thereby forcing Staller to play the vertex $v$ on her first move. All vertices are now totally dominated, except for the vertex $v$. Hence, after Dominator indicates the vertex $v$ on his second move, Staller plays any neighbor of $v$ on her second move, and the game is over, and so $\gamma_t^i(B_k) = 2$. However, the $k$ edges in the triangles that do not contain the vertex~$v$ form an induced matching in $B_k$, implying that $\indM(B_k) = k$.

\begin{prop}
\label{prop:ind-match2}
For $k \ge 1$ we have $\gamma_t^i(B_k) = 2$ and $\indM(B_k) = k$.
\end{prop}

By Proposition~\ref{prop:ind-match1}, there exist connected graphs $G$ such that $\gamma_t^i(G)$ can be arbitrarily larger than $2\indM(G)$, while by Proposition~\ref{prop:ind-match2}, there exist connected graphs $H$ such that $2\indM(H)$ can be arbitrarily larger than $\gamma_t^i(H)$. Thus, in general there is no relation between the indicated total domination number, $\gti(G)$, of an isolate-free graph $G$ and twice the induced matching number, $2\indM(G)$, of $G$. However, we pose the following question.

\begin{quest}
\label{quest:ind-math3}
Is it true that if $G$ is an isolate-free bipartite graph, then $\gti(G) \le 2\indM(G)$?
\end{quest}

For $k \ge 1$, let $J_k$ be the graph obtained from $k$ vertex disjoint copies of $C_4$ by selecting one vertex from each copy of $C_4$ and identifying these vertices into one new vertex~$v$, and then adding a new vertex $u$ and adding the edge $uv$. Let $Q_1, \ldots, Q_k$ be the $k$ copies of $C_4$ in the graph $J_k$, and so each cycle $Q_i$ contains the vertex $v$ for $i \in [k]$. Let $v_i$ be the vertex in $Q_i$ that is not adjacent to $v$ for $i \in [k]$. The graph $J_4$, for example, is illustrated in Figure~\ref{fig:B4}(b). Dominator can force the game to be completed after $k+1$ moves by indicating the vertex $v_i$ on his $i$th move for $i \in [k]$, and then indicating the vertex $u$ on his $(k+1)$st move. This strategy of Dominator forces Staller to play the vertex~$v$ and one neighbor of the vertex $v_i$ for all $i \in [k]$, thereby producing a TD-set in the graph $J_k$, implying that $\gamma_t^i(J_k) \le k+1$. (One can readily show that Staller has a strategy to prolong the game by at least~$k+1$ moves, implying that $\gamma_t^i(J_k) = k+1$.) However, selecting $k$ edges of $J_k$, one edge incident with each of the vertices $v_i$ for $i \in [k]$, produces an induced matching in $J_k$, implying that $\indM(J_k) = k$. Hence, $2\indM(J_k) - \gamma_t^i(J_k) = 2k - (k+1) = k - 1$. We state this formally as follows.

\begin{prop}
\label{prop:ind-match3}
There exist connected bipartite graphs $G$ such that $2\indM(G)$ can be arbitrarily larger than $\gamma_t^i(G)$.
\end{prop}

The following result gives a partial answer to Question~\ref{quest:ind-math3} in the case when $T$ is a tree that contains a maximum induced matching satisfying certain properties.

\begin{prop}
\label{prop:ind-math1}
If a nontrivial tree $T$ contains a maximum induced matching $M$ such that one end of every edge in $M$ is a leaf of $T$, then $\gti(T) \le 2\indM(T)$.
\end{prop}
\proof Let $T$ be a nontrivial tree of order~$n$ that contains a maximum induced matching $M$ such that one end of every edge in $M$ is a leaf of $T$. If $n = 2$, then $T = K_2$ and $\gti(T) = 2 = 2\indM(T)$. Hence, we may assume that $n \ge 3$. We define $V(M)$ as the set of vertices that are incident with an edge of $M$. Let $M = \{e_1,\ldots,e_k\}$ and let $e_i = u_iw_i$ for $i \in [k]$. By supposition, one of $u_i$ and $w_i$ is a leaf for all $i \in [k]$. Renaming vertices if necessary, we may assume that $u_i$ is a leaf for all $i \in [k]$. Let $U = \{u_1,\ldots,u_k\}$ and let $W = \{w_1,\ldots,w_k\}$, and so $V(M) = U \cup W$. Let $X$ be the boundary of the set $V(M)$ in the tree $T$, and so $X$ is the set of vertices not in the set $V(M)$ that have a neighbor in $V(M)$. We note that since every vertex in the set $U$ is a leaf (with a neighbor in $W$), the set $W$ totally dominates the boundary $X$ (and totally dominates the set $U$).  Let $Y$ be the set of vertices not totally dominated by the set $W$ in the tree $T$. We proceed further with the following claim.

\begin{claim}
\label{c:claim1}
If $Y = \emptyset$, then $\gti(T) \le 2\indM(T)$.
\end{claim}
\proof Suppose that $Y = \emptyset$, implying that $V(T) = U \cup W \cup X$ and that $W$ is a dominating set of $T$. In this case, Dominator's strategy is to indicate on his $i$th move the leaf $u_i$ for $i \in [k]$. Thus on his first $k$ moves, \D indicates the leaves $u_1,\ldots,u_k$ in turn. This forces Staller to play the unique neighbor of $u_i$, namely vertex~$w_i$, on her $i$th move for all $i \in [k]$. After these $k$ moves of Staller have been played, the vertices $w_1,\ldots,w_k$ have been added to the set that will grow to a TD-set upon completion of the game. Thus at this stage of the game, $W$ is the set of vertices that have been played. By our earlier observations, the set $W$ totally dominates the set $U \cup X$. Thus the only vertices not yet totally dominated in the game are those vertices that belong to the set $W$.

Dominator now ensures that the current set of played vertices, namely $W$, can be extended to a TD-set of the tree $T$ by adding to the set at most~$k$ vertices. This goal he readily achieves by indicating the vertices in $W$ in turn that are not yet totally dominated. More precisely, Dominator indicates the vertex $w_1$ on his $(k+1)$st move. After Staller's reply, if the game is not yet over, then Dominator indicates on his next move the vertex $w_j$ with smallest subscript $j$ that has not yet been totally dominated. Continuing in this way, Dominator has a strategy to complete the game in at most~$2k = 2\indM(T)$ moves, and so $\gti(T) \le 2k = 2\indM(T)$.~\smallqed

\medskip
By Claim~\ref{c:claim1}, we may assume that $|Y| \ge 1$, for otherwise the desired result follows. If $Y$ is not an independent set, then we could add to $M$ an arbitrary edge that belongs to the induced subtree $T[Y]$ of $T$, contradicting the maximality of the induced matching $M$. Hence, $Y$ is an independent set, and so all neighbors of vertices in $Y$ belong to the set $X$.

Dominator now employs an opening game strategy in the game that ensures that all vertices in $Y$ are totally dominated as follows. On his first more, he indicates an arbitrary vertex $y_1 \in Y$. Staller must reply by playing a neighbor of $y_1$, which as observed earlier belongs to the set $X$. Let $x_1$ be the vertex played by Staller in response to Dominator indicating the vertex $y_1$, and so $x_1 \in X$. If a vertex in $Y$ is not totally dominated by the vertex $x_1$, then Dominator indicates a vertex $y_2 \in Y$ that is not adjacent to $x_1$.  Let $x_2$ be the vertex played by Staller in response to Dominator indicating the vertex $y_2$, and so $x_2 \in X$ and $x_1 \ne x_2$. Continuing in this way, Dominator indicates on each of his next moves a vertex in $Y$ not yet totally dominated by the set of vertices played to date, thereby forcing Staller to respond by playing a vertex in $X$ that totally dominates the indicated vertex. Suppose that this process takes $r$ moves to reach the situation when all vertices of $Y$ are totally dominated. Further, suppose $y_i$ is the vertex indicated by Dominator on his $i$th move and that $x_i$ is the vertex played by Staller on her $i$th move for $i \in [r]$. We call the vertex $y_i$ the partner of the vertex $x_i$ for $i \in [r]$. This completes the opening game strategy of Dominator.

Let $X_1 = \{x_1,\ldots,x_r\} \subseteq X$ be the resulting set of vertices played by Staller on her first $r$ moves. The set $X_1$ totally dominates the set $Y$. Let $Y_1 = \{y_1,\ldots,y_r\}$ be the set of vertices indicated by Dominator on his first $r$ moves, and so $x_iy_i$ is an edge of $T$ (and $y_i$ is the partner of $x_i$) for $i \in [r]$. Let $W_1$ be the set of all vertices in $W$ that are (totally) dominated by the set $X_1$.

\begin{claim}
\label{c:claim2}
$|X_1| \le |W_1|$.
\end{claim}
\proof Suppose, to the contrary, that $|X_1| > |W_1|$.  Then there must exist $w \in W_1$ such that $|N_T(w) \cap X_1| \ge 2$. Renaming vertices if necessary, we may assume that $x_1,x_2 \in N_T(w) \cap X_1$. But then, the induced matching obtained from $M$ by removing the edge $uw$, where $u \in U$ is the neighbor of $w$, and adding the edges $x_1y_1$ and $x_2y_2$ has cardinality $|M| + 1$, a contradiction to the maximality of $M$.~\smallqed

\medskip
By Claim~\ref{c:claim2}, $|X_1| \le |W_1|$. Let $X_2 = X \setminus X_1$ and let $W_2 = W \setminus W_1$. Further, let $U_i$ be the set of neighbors of vertices in $W_i$ that belong to the set $U$ for $i \in [2]$. Hence, the induced subgraph $T[U_i \cup W_i]$ is a disjoint union of $|W_i|$ copies of $K_2$, where each copy of $K_2$ consists of a vertex in $W_i$ and its unique neighbor that belongs to the set $U$ for $i \in [2]$.

In his middle game strategy, Dominator ensures that all vertices in $W$ are played as follows. Dominator indicate on his $(r+i)$th move the leaf $u_i$ for $i \in [k]$. Thus on his next $k$ moves immediately following the first $r$ moves in his opening game strategy, \D indicates the leaves $u_1,\ldots,u_k$ in turn. This forces Staller to play the unique neighbor of $u_i$, namely vertex~$w_i$, on her $(r+i)$th move for all $i \in [k]$. After these $r+k$ moves of Staller have been played, the vertices $w_1,\ldots,w_k$ have been added to the set of played vertices, yielding the current set $X_1 \cup W$  of vertices that have been played to date by Staller. The set $W$ of vertices played to date in the game totally dominates the set $U \cup X$. Recall that Dominator's opening game strategy ensures that the set $X_1$ totally dominates the set $Y \cup W_1$. Thus, Dominator middle game strategy adds $k$ new vertices to the set of vertices played by Staller, in addition to the $r$ vertices she played in the opening game phase of the game. Upon completion of his middle game strategy, the only vertices not yet totally dominated in the game are those vertices that belong to the set $W_2$.

In his end game strategy, Dominator ensures that all vertices in $W_2$ are totally dominated as follows. Dominator orders the vertices in $W_2$ sequentially. Let $w_{2,1},w_{2,2},\ldots,w_{2,q}$ be the resulting ordering of the vertices in $W_2$ by Dominator, where $q = |W_2|$. On each of his subsequent moves, Dominator  indicates the vertex $w_{2,j}$ with smallest subscript $j$ that has not yet been totally dominated where $j \in [q]$. Thus, Dominator's end game strategy adds $q$ new vertices to the set of vertices played by Staller, in addition to the $r$ vertices she played in the opening game phase of the game and the $k$ vertices she played in the middle game phase of the game. Upon completion of Dominator's end game strategy, at most $r + k + q$ vertices are played by Staller, and the resulting set of played vertices is a TD-set of the tree $T$. Recall that $|W_1| \ge |X_1| = r$, and so $q = |W_2| = |W| - |W_1| \le k - r$. Therefore, Dominator has a strategy to complete the game in at most~$(r + q) + k \le 2k = 2|M| = 2\indM(T)$ moves, and so $\gti(T) \le 2k = 2\indM(T)$.~\QED

\section{The game total domination number}

We show in this section that in general there is no relation between the indicated total domination number, $\gti(G)$, of an isolate-free graph $G$ and the game total domination number, $\gtd(G)$, of $G$.

\begin{figure}[ht!]
\begin{center}
\begin{tikzpicture}[scale=0.7,style=thick,x=0.9cm,y=0.9cm]
\def\vr{3pt} 
\path (-1,0) coordinate (a);
\path (1,2) coordinate (x1); \path (1,1) coordinate (x2); \path (1,-2) coordinate (x3);
\path (3,2) coordinate (y1); \path (3,1) coordinate (y2); \path (3,-2) coordinate (y3);


\foreach \i in {1,2,3}
{\draw (x\i)--(y\i); \draw (a)--(x\i); }

\draw (a) [fill=white] circle (\vr);

\foreach \i in {1,2,3}
{  \draw (x\i)  [fill=white] circle (\vr); \draw (y\i)  [fill=white] circle (\vr); }

\draw (1,.0) [fill=black] circle (1.5pt);
\draw (1,-.5) [fill=black] circle (1.5pt);
\draw (1,-1) [fill=black] circle (1.5pt);
\draw (3,0) [fill=black] circle (1.5pt);
\draw (3,-.5) [fill=black] circle (1.5pt);
\draw (3,-1) [fill=black] circle (1.5pt);

\draw[anchor = east] (a) node {{\small $v$}};
\draw[anchor = south] (x1) node {{\small $x_1$}};
\draw[anchor = south] (x2) node {{\small $x_2$}};
\draw[anchor = north] (x3) node {{\small $x_k$}};
\draw[anchor = south] (y1) node {{\small $y_1$}};
\draw[anchor = south] (y2) node {{\small $y_2$}};
\draw[anchor = north] (y3) node {{\small $y_k$}};

\end{tikzpicture}
\caption{The tree $T_k$ in the proof of Proposition~\ref{prop:subdivide-star}}
\label{fig:SubdividedStar}
\end{center}
\end{figure}
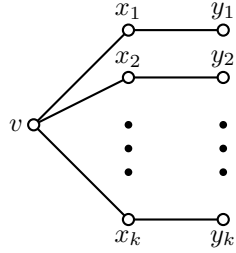

\begin{prop}
\label{prop:subdivide-star}
If $k \ge 3$ and $T_k$ is the tree obtained from a star $K_{1,k}$ by subdividing every edge once, then $\gamma_t^i(T_k) = \Gamma_t(T_k) = 2k$ and $\gtd(T_k)=k+1$.
\end{prop}
\proof For $k \ge 3$, let $T = T_k$ be the tree obtained from a star $K_{1,k}$ by subdividing every edge exactly once. Let $v$ be the central vertex of $T$ (of degree~$k)$, and let $N_T(v) = \{x_1,x_2,\ldots,x_k\}$. Further, let $y_i$ be the leaf neighbor of $x_i$ for $i \in [k]$. The tree $T_k$ is illustrated in Figure~\ref{fig:SubdividedStar}.

The set of leaves together with the set of support vertices in $T_k$ is a minimal TD-set of maximum cardinality, and thus $\Gamma_t(T_k) = 2k$.

When the total domination game is played on $T_k$, Dominator plays the vertex $v$ as his first move.  The set of legal moves remaining is the set of support vertices, namely $\{x_1,x_2,\ldots,x_k\}$ and each of these vertices must be played during the game in order to totally dominate the leaves of $T_k$. This strategy of Dominator shows that $\gtd(T_k) \le k+1$. In any play of the total domination game on $T_k$ exactly $k$ moves will be required to totally dominate the set of leaves, namely the set $\{y_1,\ldots, y_k\}$. Furthermore, either $v$ or $y_1$ must be played to totally dominate the vertex $x_1$.  This shows that $\gtd(T_k) \ge k+1$. Therefore, $\gtd(T_k) = k+1$.

By Proposition~\ref{prop:upper-Tdom} we have $\gamma_t^i(T_k) \ge \Gamma_t(T_k)=2k$. By Proposition~\ref{prop:ind-math1} we readily infer that $\gamma_t^i(T_k) \le 2k$. Consequently, $\gamma_t^i(T_k)=2k$.~\QED

\medskip
By Proposition~\ref{prop:subdivide-star}, there exist connected graphs $G$ such that $\gamma_t^i(G)$ can be arbitrarily larger than $\gtd(G)$.  Using the result of Proposition~\ref{prop:subdivide-star} and the fact that $\Gamma_t(G) \le \OOIR(G)$, we see that there exist trees $T$ such that $\OOIR(T)$ is arbitrarily larger than $\gtd(T)$. We show next that there exist connected graphs $G$ such that $\gtd(G)$ can be arbitrarily larger than both $\gamma_t^i(G)$ and $\OOIR(G)$.

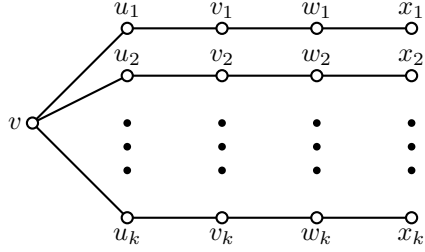
\begin{figure}[ht!]
\begin{center}
\begin{tikzpicture}[scale=0.7,style=thick,x=0.9cm,y=0.9cm]
\def\vr{3pt} 
\path (-1,0) coordinate (a);
\path (1,2) coordinate (x1);
\path (1,1) coordinate (x2);
\path (1,-2) coordinate (x3);
\path (3,2) coordinate (y1);
\path (3,1) coordinate (y2);
\path (3,-2) coordinate (y3);
\path (5,2) coordinate (z1);
\path (5,1) coordinate (z2);
\path (5,-2) coordinate (z3);
\path (7,2) coordinate (w1);
\path (7,1) coordinate (w2);
\path (7,-2) coordinate (w3);


\foreach \i in {1,2,3}
{\draw (x\i)--(y\i); \draw (a)--(x\i); \draw (y\i)--(z\i); \draw (z\i)--(w\i);}

\draw (a) [fill=white] circle (\vr);

\foreach \i in {1,2,3}
{  \draw (x\i)  [fill=white] circle (\vr); \draw (y\i)  [fill=white] circle (\vr);
\draw (z\i)  [fill=white] circle (\vr); \draw (w\i)  [fill=white] circle (\vr); }

\draw (1,0) [fill=black] circle (1.5pt);
\draw (1,-.5) [fill=black] circle (1.5pt);
\draw (1,-1) [fill=black] circle (1.5pt);
\draw (3,0) [fill=black] circle (1.5pt);
\draw (3,-.5) [fill=black] circle (1.5pt);
\draw (3,-1) [fill=black] circle (1.5pt);
\draw (5,0) [fill=black] circle (1.5pt);
\draw (5,-.5) [fill=black] circle (1.5pt);
\draw (5,-1) [fill=black] circle (1.5pt);
\draw (7,0) [fill=black] circle (1.5pt);
\draw (7,-.5) [fill=black] circle (1.5pt);
\draw (7,-1) [fill=black] circle (1.5pt);

\draw[anchor = east] (a) node {{\small $v$}};
\draw[anchor = south] (x1) node {{\small $u_1$}};
\draw[anchor = south] (x2) node {{\small $u_2$}};
\draw[anchor = north] (x3) node {{\small $u_k$}};
\draw[anchor = south] (y1) node {{\small $v_1$}};
\draw[anchor = south] (y2) node {{\small $v_2$}};
\draw[anchor = north] (y3) node {{\small $v_k$}};
\draw[anchor = south] (z1) node {{\small $w_1$}};
\draw[anchor = south] (z2) node {{\small $w_2$}};
\draw[anchor = north] (z3) node {{\small $w_k$}};
\draw[anchor = south] (w1) node {{\small $x_1$}};
\draw[anchor = south] (w2) node {{\small $x_2$}};
\draw[anchor = north] (w3) node {{\small $x_k$}};

\end{tikzpicture}
\caption{The tree $T_k$ in the proof of Proposition~\ref{prop:subdivide-star3}}
\label{fig:3SubdividedStar}
\end{center}
\end{figure}

\begin{prop}
\label{prop:subdivide-star3}
For $k \ge 4$ an even integer, if $T_k$ is the tree obtained from a star $K_{1,k}$ by subdividing every edge three times, then $\gamma_t^i(T_k)\le 2k+2 = \OOIR(T_k)$ and $\gtd(T_k) \ge \frac{5}{2}k$.
\end{prop}
\proof For $k \ge 4$ an even integer, let $T_k$ be the tree obtained from a star $K_{1,k}$ by subdividing every edge exactly three times. Let $v$ be the central vertex of $T_k$ (of degree~$k)$, and let $Q_i \colon vu_iv_iw_ix_i$ be the $k$ paths of length~$4$ emanating from $v$ in $T_k$ for $i \in [k]$. The tree $T_k$ is illustrated in Figure~\ref{fig:3SubdividedStar}.

First we show that $\OOIR(T_k) = 2(k + 1)$. Let $M$ be an induced matching in $T_k$ of (maximum) cardinality~$\indM(T_k)$. By the maximality of $M$, the induced matching $M$ contains an edge from the set $\{u_iv_i,v_iw_i,w_ix_i\}$ for all $i \in [k]$. If the edge $u_iv_i$ or the edge $v_iw_i$ belong to the induced matching $M$, then we can replace such an edge in $M$ with the edge $w_ix_i$. Hence, we may choose the induced matching $M$ to contain the edge $w_ix_i$ for all $i \in [k]$. No additional edge incident with $v_i$ or $w_i$ can be added to these $k$ edges without violating the requirement that $M$ is an induced matching. Hence, $M$ contains exactly one additional edge, namely an edge incident with the vertex~$v$. Hence, the set
\[
M = \{vu_1\} \cup \bigcup_{i=1}^k\{w_ix_i\},
\]
for example, is an induced matching in $T_k$ of maximum cardinality, and so $\indM(T_k) = |M| = k + 1$. By Theorem~\ref{thm:cube-paper}, we therefore infer that $\OOIR(T_k) = 2(k + 1)$.

We show next that $\gtd(T_k) \ge \frac{5}{2}k$. Let the total domination D-game be played on $T_k$.  We provide a strategy for Staller that forces at least $\frac{5}{2}k$ moves to be made in the game. We note that when the game has ended, the vertex $w_i$ has been played in order to totally dominate the leaf $x_i$ for each $i \in [k]$. Further, at least one vertex from $\{v_i,x_i\}$ will have been played in order to totally dominate the vertex $w_i$ for each $i \in [k]$. Thus at least two vertices are played from every set $\{v_i,w_i,x_i\}$ for all $i \in [k]$. Staller's strategy is to play vertex $u_i$, for as many values of $i \in [k]$ as possible, whenever no other vertex from $\{u_i,v_i,w_i,x_i\}$ has been played. Since $k$ is even, by following this strategy she can ensure that at least $k/2$ of the vertices in $\{u_1,\ldots,u_k\}$ will be played.   For each $j \in [k]$ for which Staller played vertex $u_j$, as observed earlier the vertex $w_j$ and at least one of $v_j$ and $x_j$ will be played in the course of the game, implying that at least three vertices from the set $\{u_j,v_j,w_j,x_j\}$ are played. Therefore upon completion of the total domination D-game, Staller guarantees that at least three moves are played from at least~$k/2$ of the sets $\{u_j,v_j,w_j,x_j\}$, and from the remainder of the sets $\{u_j,v_j,w_j,x_j\}$ at least two moves are played. Therefore, when the game is complete at least $3 \times \frac{k}{2} + 2 \times \frac{k}{2} = \frac{5}{2}k$ moves were made, implying that $\gtd(T_k) \ge \frac{5}{2}k$.

Finally, let the indicated total domination game be played on $T_k$.  We give a strategy for Dominator that ensures at most $2k+2$ vertices are chosen. In his first $k$ moves Dominator indicates - in any order - the vertices in $\{x_1,\ldots,x_k\}$, thereby forcing \St to respond by playing all $k$ support vertices, namely the vertices in the set $\{w_1,\ldots,w_k\}$. After Staller's first $k$ moves the result is that the set $\cup_{i=1}^k\{v_i,x_i\}$ is therefore totally dominated.  Now Dominator indicates (in order) the vertices in the sequence $u_1, \ldots,u_k$.  If Staller responds by playing $v_1, \ldots, v_k$, respectively, then every vertex except for vertex $v$ has been totally dominated. As a result $2k+1$ vertices will be chosen by Staller. On the other hand, if Staller plays $v_1,\ldots, v_j$, respectively, for some $1 \le j<k$ and then plays vertex~$v$ when \D indicates the vertex $u_{j+1}$, exactly one vertex from $\{v_i,x_i\}$ will be played by Staller for each $i \in [k] \setminus [j]$ when the game ends.  In this case $k + j + 1 + (k-j) = 2k+1$ vertices will be chosen upon completion of the game.  Finally, if Staller plays the vertex $v$ when Dominator indicates vertex $u_1$, then exactly one vertex from $\{v_i,x_i\}$ will be played by Staller (in order to totally dominate the vertex $w_i$) for each $i \in [k]$ and one vertex from  $\{u_1,\ldots, u_k\}$  will be played by Staller
when Dominator indicates the vertex $v$. In this case a total of $k+1+k+1=2k+2$ vertices will be chosen by Staller. This strategy by Dominator ensures that at most $2k+2$ vertices will be played by Staller. Therefore,  $\gamma_t^i(T_k) \le 2k+2$.~\QED

\medskip
By Proposition~\ref{prop:subdivide-star3}, there therefore exist connected graphs $G$ such that $\gtd(G)$ can be arbitrarily larger than $\gamma_t^i(G)$.

We also remark that when the total domination game (D-game) is played on the graph $F_k$ defined earlier immediately before the statement of Proposition~\ref{prop:cartesian} where $k \ge 5$, then Dominator can ensure that both vertices $u_1$ and $v_1$ are played in the first three moves. Indeed, \D plays the vertex $u_1$ in his first move. Staller can force the game to last at least three moves by playing vertex $u_k$ on her first move. Therefore, $\gtd(F_k) = 3$. Recall that by Proposition~\ref{prop:cartesian}, we have $\gamma_t^i(F_k)=4$. Hence for $k \ge 5$, the graph $F_k$ is another example of a graph satisfying $\gtd(F_k) < \gamma_t^i(F_k)$.

\section{The Grundy total domination number}

Recall that the length of a longest total dominating sequence in $G$ is the Grundy total domination number, $\grt(G)$, of~$G$. The definition of the indicated total domination game implies that the sequence of vertices selected by Staller in the indicated total domination game is a total dominating sequence. From this we infer that $\gti(G) \le \grt(G)$. Let $A$ be an open-open irredundant set in $G$ of cardinality $\OOIR(G)$ and construct a sequence $S$ using the vertices of $A$ in any order.  If $S$ is not a total dominating sequence, then it can be extended to form one.  This shows that $\OOIR(G) \le \grt(G)$. We state this formally as follows.

\begin{observation}
\label{ob:Grundy-Tdom}
If $G$ is an isolate-free graph, then $\gti(G) \le \grt(G)$ and $\OOIR(G) \le \grt(G)$.
\end{observation}

As a consequence of Proposition~\ref{prop:path} and the following result given in~\cite{BrHeRa-16}, we infer that there exist trees $T$ such that $\grt(T)$ is arbitrarily larger than $\gamma_t^i(T)$.

\begin{theorem}{\rm (\cite{BrHeRa-16})}
\label{thm:Grundy}
If $T$ is a tree of order~$n$, then $\grt(T) = n$ if and only if $T$ has a perfect matching.
\label{cor:tree}
\end{theorem}

\section{Summary and concluding remarks}

\begin{figure}[ht!]
	\begin{center}
		\begin{tikzpicture}
			\node (2) at (1.5,3) {$2\indM$};
			\node (3) at (0,3) {$\gamma_t$};
			\node (4) at (0,4) {$\Gamma_t$};
			\node (5) at (-1.5,5) {$\gti$}; 
			\node (6) at (0,5) {$\OOIR$};
			\node (7) at (1.5,5) {$\gtd$};
			\node (8) at (0,6) {$\grt$}; 			
			\draw (3) -- (4) -- (5) -- (8);
			\draw (4) -- (6) -- (8);
			\draw (3) -- (7) -- (8);
			\draw (2) -- (6);
		\end{tikzpicture}
	\caption{Relations between the invariants studied in this paper}
	\label{fig:hasse}
	\end{center}
\end{figure}
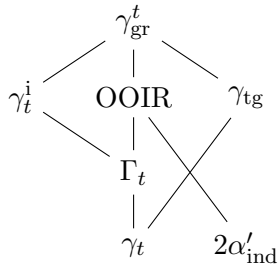

The Hasse diagram in Figure~\ref{fig:hasse} shows the relationships between the invariants studied in this paper for graphs with no isolated vertices. The invariants $\gti$, $\OOIR$ and $\gtd$ are pairwise incomparable.  Indeed, Propositions~\ref{prop:join} and~\ref{prop:cartesian} show that the differences
$\gti(G)-\OOIR(G)$ and  $\OOIR(G)-\gti(G)$ can be arbitrarily large; Propositions~\ref{prop:subdivide-star} and~\ref{prop:subdivide-star3} show that $\gti(G)-\gtd(G)$ and $\gtd(G)-\gti(G)$ can be arbitrarily large;  Propositions~\ref{prop:subdivide-star} and~\ref{prop:subdivide-star3} together with the fact that every minimal TD-set is open-open irredundant show that $\OOIR(G)-\gtd(G)$ and $\gtd(G)-\OOIR(G)$ can be arbitrarily large.  Consider the graph $B_k$ defined in Section~\ref{S:ind-math} for each positive integer $k$.  We see from Proposition~\ref{prop:ind-match2} that $\gti(B_k)=2$ and $2\indM(B_k)=2k$.   It is also easy to verify that $\gtd(B_k)=\gamma_t(B_k)=\Gamma_t(B_k)=2$.  Thus,
$2\indM$ can be arbitrarily larger than $\gti, \gtd, \gamma_t$ and $\Gamma_t$.  On the other hand, if $G$ is the corona of a complete graph of order $k$, then
$2\indM(G)=2$, $\gti(G)=\gamma_t(G)=\Gamma_t(G)=k$ and $\gtd(G)=k+1$, which shows that $2\indM$ can be arbitrarily smaller than these four invariants as well.  Therefore,
$2\indM$ is incomparable with each of $\gti, \gtd, \gamma_t$ and $\Gamma_t$.

\end{document}